\documentclass[a4paper,11pt]{article}
\usepackage{amsmath,amssymb,amsthm,a4wide}
\usepackage[utf8]{inputenc}
\usepackage[T1]{fontenc}
\usepackage{graphicx}        
\usepackage{url}
\usepackage{color}
\usepackage{caption}
\usepackage{mathrsfs}
\usepackage{algorithm,algpseudocode}
\usepackage{authblk}

\newtheorem{lemma}      {Lemma}      [section]

\newcommand{\nc}{\newcommand}

\nc{\dsp}{\displaystyle}
\nc{\mrm}{\mathrm}

\nc{\res}{\boldsymbol{r}}
\nc{\bs}{\boldsymbol{s}}
\nc{\bz}{\boldsymbol{z}}
\nc{\bw}{\boldsymbol{w}}
\nc{\br}{\boldsymbol{r}}
\nc{\bx}{\boldsymbol{x}}
\nc{\by}{\boldsymbol{y}}
\nc{\ba}{\boldsymbol{a}}
\nc{\bn}{\boldsymbol{n}}
\nc{\bg}{\boldsymbol{g}}
\nc{\bh}{\boldsymbol{h}}
\nc{\bu}{\boldsymbol{u}}
\nc{\bv}{\boldsymbol{v}}
\nc{\be}{\boldsymbol{e}}
\nc{\bE}{\boldsymbol{E}}
\nc{\bH}{\boldsymbol{H}}
\nc{\bJ}{\boldsymbol{J}}
\nc{\bbeta}{\boldsymbol{\beta}}
\nc{\bxi}{\boldsymbol{\xi}}

\nc{\curl}{\boldsymbol{curl}}
\nc{\bfX}{\mathbf{X}}
\nc{\bfU}{\mathbf{U}}
\nc{\bfV}{\mathbf{V}}
\nc{\R}{\mathbb{R}}
\nc{\mbP}{\mathbb{P}}
\nc{\mbH}{\mathbb{H}}
\nc{\mbA}{\mathbb{A}}
\nc{\mT}{\mathrm{T}}
\nc{\mP}{\mathrm{P}}
\nc{\mQ}{\mathrm{Q}}
\nc{\mR}{\mathrm{R}}
\nc{\mD}{\mathrm{D}}
\nc{\mB}{\mathrm{B}}
\nc{\mX}{\mathrm{X}}
\nc{\mA}{\mathrm{A}}
\nc{\mL}{\mathrm{L}}
\nc{\mS}{\mathrm{S}}
\nc{\mY}{\mathrm{Y}}
\nc{\mV}{\mathrm{V}}
\nc{\mH}{\mathrm{H}}
\nc{\mC}{\mathrm{C}}
\nc{\mJ}{\mathrm{J}}
\nc{\Id}{\mrm{Id}}

\nc{\mId}{\mathrm{Id}}

\nc{\lbr}{\lbrack}
\nc{\rbr}{\rbrack}
\nc{\RR}{\mathbb{R}}
\nc{\CC}{\mathbb{C}}

\nc{\ctru}{\mathfrak{u}}
\nc{\ctrv}{\mathfrak{v}}

\nc{\dof}{\mrm{dof}}

\title{\textbf{A new perspective on  the fictitious space lemma}}
\date{}
\author[1]{X.Claeys}
\affil[1]{Sorbonne Université, Université Paris-Diderot SPC, CNRS, 
  Laboratoire Jacques-Louis Lions}

\begin{document}

\maketitle

\section*{Introduction}

In the present contribution we propose a new proof of the
so-called fictitious space lemma. For the proof, we exhibit
an explicit expression for the inverse of additive Schwarz preconditionners
in terms of Moore-Penrose pseudo inverse of the map associated to
the decomposition over the subdomain partition.

We will first briefly recall the definition of the pseudo-inverse of a matrix
and some of its remarkable properties. We will then  explain how this concept
can be used to reformulate the fictitious space lemma in a very compact form.
We will then give an aletrnative proof of the fictious space lemma.
As a remarkable feature, this proof does not rely on Cauchy-Schwarz inequality,
as opposed to previous proofs provided by Nepomnyaschikh \cite{MR1189564,MR1126677},
Griebel and Oswald \cite{MR1324736} or Dolean, Jolivet and Nataf \cite{MR3450068},
see also \cite{MR1193013}. The present proof applies directly in the infinite dimensional case.

\section{Moore-Penrose pseudo-inverse}

Assume given Hilbert spaces $\mH$ (resp.$\mV$) equipped with the norms
$\Vert \cdot\Vert_{\mH}$ (resp. $\Vert \cdot\Vert_{\mV}$) and consider
a surjective map $\mR:\mV\to \mH$. Define $\mR^{-1}(\{\by\}):=\{\bx'\in\mH,\;
\mR\bx' = \by\}$. The Moore-Penrose pseudoinverse
of this map, denoted $\mR^{\dagger}:\mH\to \mV$ is defined, for all $\by\in \mH$, by
\begin{equation}\label{DefPseudoInverse}
  \mR\cdot\mR^{\dagger}\by = \by\quad\;\;\text{and}\;\;\quad
  \Vert \mR^{\dagger}\by\Vert_{\mV} = \inf_{\bx\in \mR^{-1}(\{\by\})}\Vert \bx\Vert_{\mV}.
\end{equation}
The property above directly implies that $\mR^{\dagger}$ is injective.  
Let us denote $\mV_{\mR} = \mrm{Ker}(\mR)^{\perp}$. For any $\by\in \mH$, since 
the restricted operator $\mR\vert_{\mV_{\mR}}:\mV_{\mR}\to \mH$ is a bijection,
there exists a unique $\bx\in \mV_{\mR}$ such that $\mR\bx = \by$.
Besides, if $\bx'\in \mV$ is another element satisfying $\mR\bx' = \by$
then $\bx-\bx'\in \mrm{Ker}(\mR)$ so that $(\bx,\bx-\bx')_{\mV}$ and thus,
by Pythagore's rule,
\begin{equation}
  \Vert \bx\Vert_{\mV}^{2}\leq  \Vert \bx\Vert_{\mV}^{2} + \Vert \bx-\bx'\Vert_{\mV}^{2} = \Vert \bx'\Vert_{\mV}^{2}
\end{equation}
As a consequence $\bx\in \mV$ solves the minimization
problem (\ref{DefPseudoInverse}) i.e. $\bx = \mR^{\dagger}\by$.
From this discussion we conclude that  $\mR^{\dagger} = (\mR\vert_{\mV_{\mR}})^{-1}$.

\quad\\
The property $\mR\mR^{\dagger} = \mrm{Id}$ implies that $(\mR^{\dagger}\mR)^{2} =
\mR^{\dagger}(\mR\mR^{\dagger})\mR = \mR^{\dagger}\mR$ i.e. $\mR^{\dagger}\mR$
is a projector. Because $\mR^{\dagger}$ is injective we obtain
$\mrm{Ker}(\mR^{\dagger}\mR) = \mrm{Ker}(\mR)$.Besides for any $\bx\in \mV_{R}$
satisfying $\mR\bx = \by$, we have seen that $\bx = \mR^{\dagger}\by = \mR^{\dagger}\mR\bx$
which implies $\mV_{\mR} = \mrm{Im}(\mR^{\dagger}\mR)$. As a conclusion, since
$\mrm{Ker}(\mR)$ and $\mV_{\mR}$ are orthogonal by definition, we conclude that 
$\mR^{\dagger}\mR$ is an orthogonal projection,
which rewrites
\begin{equation}\label{SelfAdjointness}
(\mR^{\dagger}\mR\bx,\by)_{\mV} = (\bx,\mR^{\dagger}\mR\by)_{\mV}\quad \forall \bx,\by\in \mV.
\end{equation}

\section{Weighted pseudo-inverse}\label{secWeightedPseudoInverse}

Keeping the notations from the previous section,
consider  continuous operator $\mB:\mV\to \mV$, and
assume this operator is self-adjoint so that
it induces a scalar product $(\bx,\by)_{\mB}:=(\mB\bx,\by)_{\mV}$
and a norm $\Vert \bx\Vert_{\mB}:=\sqrt{(\bx,\bx)_{\mB}}$. 
To each such $\mB$ can be associated a so-called "weighted pseudo-inverse"
$\mR^{\dagger}_{\mB}:\mH\to \mV$ defined, for all $\forall\by\in \mH$  by 
\begin{equation}\label{WeightedPseudoInverse}
  \mR\cdot\mR^{\dagger}_{\mB}\by = \by\quad\;\;\text{and}\quad\;\;
  \Vert \mR^{\dagger}_{\mB}\by\Vert_{\mB} = \inf_{\bx\in \mR^{-1}(\{\by\})}\Vert \bx\Vert_{\mB}.
\end{equation}
The operator $\mR^{\dagger}_{\mB}$ satisfies the same
properties as $\mR^{\dagger}$ except that $(\;,\;)_{\mV}$ is this
times replaced by $(\;,\;)_{\mB}$. In particular (\ref{SelfAdjointness}) rewrites
$(\mR^{\dagger}_{\mB}\mR\bx,\by)_{\mB} = (\bx,\mR^{\dagger}_{\mB}\mR\by)_{\mB}$
for all $\bx,\by\in \mV$. Taking account of the expression of $(\;,\;)_{\mB}$
this is equivalent to
\begin{equation}\label{WeightedAdjointness}
  \mB\,\mR^{\dagger}_{\mB}\mR = (\mR^{\dagger}_{\mB}\mR)^{*}\mB
\end{equation}
where, for any continuous linear operator $\mrm{M}:\mV\to \mV$ we denote $\mrm{M}^{*}$
its adjoint with respect to $(\;,\;)_{\mV}$ defined by
$(\mrm{M}\bx,\by)_{\mV} = (\bx,\mrm{M}^{*}\by)_{\mV}$ for all $\bx,\by\in\mV$.
Property (\ref{WeightedAdjointness}) leads to a lemma.

\begin{lemma}\label{FormuleInverse}\quad\\
  $\mR\mB^{-1}\mR^{*} = (\,(\mR_{\mB}^{\dagger})^{*}\mB\mR_{\mB}^{\dagger}\,)^{-1}$
\end{lemma}
\noindent \textbf{Proof:}

Since $\mR\mR^{\dagger}_{\mB} = \mrm{Id}$ by construction, the lemma is
a consequence of (\ref{WeightedAdjointness}) through  direct calculation
$  (\mR\mB^{-1}\mR^{*})\cdot((\mR_{\mB}^{\dagger})^{*}\mB\mR_{\mB}^{\dagger}) 
= \mR\mB^{-1}(\mR_{\mB}^{\dagger}\mR)^{*}\mB\mR_{\mB}^{\dagger}
= \mR(\mB^{-1}\mB)\mR_{\mB}^{\dagger}\mR\mR_{\mB}^{\dagger} = (\mR\mR_{\mB}^{\dagger})^{2} = \mrm{Id}$.
\hfill$\Box$

\section{Re-interpretation of the fictious space lemma}\label{ProofFictitiousSpaceLemma}

In this section, we provide a new proof of the fictitious space lemma relying
on the concept weighted pseudoinverse. As a preliminary, let us recall a classical
caracterisation of extremal eigenvalues of self-adjoint operators
(see e.g. theorem 1.2.1 and theorem 1.2.3  in \cite{MR1417493}).

\begin{lemma}\label{Minimax}\quad\\
  Assume $\mH$ is an Hilbert space equipped with the
  scalar product  $(\;,\;)_{\mH}$ and let 
  $\mT:\mH\to \mH$ be a bounded operator that is self-adjoint for
  $(\;,\;)_{\mH}$. Denoting $\sigma(\mT)$ the spectrum of
  $\mT$, we have
  \begin{equation*}
    \inf \sigma(\mT) = \inf_{\bx\in \mH\setminus\{0\}}\frac{(\mT\bx,\bx)_{\mH}}{(\bx,\bx)_{\mH}}
    \quad\quad\quad
    \sup \sigma(\mT) = \sup_{\bx\in \mH\setminus\{0\}}\frac{(\mT\bx,\bx)_{\mH}}{(\bx,\bx)_{\mH}}    
  \end{equation*}
\end{lemma}

\noindent 
This lemma holds independently of the choice of the scalar product,
provided that $\mT$ be self-adjoint with respect to it. As a consequence of
the previous lemma, if  $\alpha(\;,\;)$
and $\beta(\;,\;)$ are two scalar products over $\mH$ and $\mT$
is self-adjoint with respect to both, then
$\inf_{\bx\in \mH\setminus\{0\}}\alpha(\mT\bx,\bx)/\alpha(\bx,\bx) =
\inf_{\bx\in \mH\setminus\{0\}}\beta(\mT\bx,\bx)/\beta(\bx,\bx)$, and a similar
result holds for the supremum.

\quad\\
Now we recall the fictitious space lemma, adopting the same formulation
of this result as \cite[Lemma 7.4]{MR3450068} and \cite[p.168]{MR1324736}.

\begin{lemma}\label{FictitiousSpaceLemma}\quad\\
  Let $\mH$ and $\mV$ be two Hilbert spaces equipped with the
  scalar products $(\;,\;)_{\mH}$ and  $(\;,\;)_{\mV}$.
  Let $\mA:\mH\to \mH$ (resp. $\mB:\mV\to \mV$) be a
  bounded operator that is positive definite self-adjoint
  with respect to $(\;,\;)_{\mH}$ (resp. $(\;,\;)_{\mV}$), and denote
  $(\bu,\bv)_{\mA}:=(\mA\bu,\bv)_{\mH}$ (resp. $(\bu,\bv)_{\mB}:=(\mB\bu,\bv)_{\mV}$).
  Suppose that there exists a surjective bounded linear operator $\mR:\mV\to \mH$,
  and constants $c_{\pm}>0$ such that
  \begin{itemize}
  \item[i)] for all $\bu\in \mH$ there exists $\bv\in \mV$
    with $\mR\bv = \bu$ and $c_{-}(\bv,\bv)_{\mB}\leq (\bu,\bu)_{\mA}$,  
  \item[ii)] $(\mR\bv,\mR\bv)_{\mA}\leq c_{+}(\bv,\bv)_{\mB}$ for all $\bv\in \mV$.\\[-5pt]
  \end{itemize}
  Then, denoting $\mR^{*}:\mH\to \mV$ the linear map defined by $(\mR\bu,\bv)_{\mH} = (\bu,\mR^{*}\bv)$
  for all $\bu\in\mV,\bv\in\mH$, we have
  \begin{equation}\label{SpectralEstimate}
    c_{-}(\bu,\bu)_{\mA}\leq (\mR\mB^{-1}\mR^{*}\mA\bu,\bu)_{\mA}
    \leq c_{+}(\bu,\bu)_{\mA}\quad \forall \bu\in\mH.
  \end{equation}
  In addition, if $c_{\pm}$ are the optimal constants satisfying \textit{i)-ii)}
  then the bounds in (\ref{SpectralEstimate}) are optimal as well.
\end{lemma}

\noindent 
\textbf{Proof:}

We simply reformulate \textit{i)-ii)} by means of the weighted pseudo-inverse. If
\textit{i)} holds then, for any $\bu\in\mH$ we have
$c_{-}\Vert\bv\Vert_{\mB}^{2} \leq (\bu,\bu)_{\mA}\forall \bv\in \mR^{-1}(\{\bu\})$.
Taking the inifimum and using (\ref{WeightedPseudoInverse}), we obtain $c_{-}\Vert\mR_{\mB}^{\dagger}\bu\Vert_{\mB}
\leq (\bu,\bu)_{\mA}$. On the other hand, it is clear that, if
$c_{-}\Vert\mR_{\mB}^{\dagger}\bu\Vert_{\mB}\leq (\bu,\bu)_{\mA}\forall\bu\in \mH$ then
\textit{i)} holds.

Next if \textit{ii)} holds, then we have $(\bu,\bu)_{\mA}\leq c_{+}\Vert\bv\Vert_{\mB}\;
\forall \bv\in \mR^{-1}(\{\bu\})$ and for all $\bu\in \mH$. Taking the infinimum over
$ \bv\in \mR^{-1}(\{\bu\})$ and using (\ref{WeightedPseudoInverse}), we conclude that 
$(\bu,\bu)_{\mA}\leq c_{+}\Vert\mR_{\mB}^{\dagger}\bu\Vert_{\mB}^{2}\;\forall \bu\in \mH$,
and this is equivalent due to the optimality condition in (\ref{WeightedPseudoInverse}).
To conclude we have just shown that conditions \textit{i)-ii)} in Lemma \ref{FictitiousSpaceLemma}
are actually equivalent to
\begin{equation}\label{SpectralEstimat2}
  c_{-}(\mR_{\mB}^{\dagger}\bu,\mR_{\mB}^{\dagger}\bu)_{\mB}\leq (\bu,\bu)_{\mA}\leq
  c_{+}(\mR_{\mB}^{\dagger}\bu,\mR_{\mB}^{\dagger}\bu)_{\mB}\quad \forall\bu\in\mH.
\end{equation}
Next define $\mS:=(\mR_{\mB}^{\dagger})^{*}\mB\mR_{\mB}^{\dagger}$, which is obviously bounded positive
definite self-adjoint so it induces a scalar
product $(\bu,\bv)_{\mS}:= (\mS\bu,\bv)_{\mH}$ and a norm $\Vert \bu\Vert_{\mS}:=\sqrt{(\bu,\bu)_{\mS}}$.
We can re-write $(\mR_{\mB}^{\dagger}\bu,\mR_{\mB}^{\dagger}\bu)_{\mB} = (\bu,\bu)_{\mS}$, and
$(\bu,\bu)_{\mA} = (\mR\mB^{-1}\mR^{*}\mA\bu,\bu)_{\mS}$ according to Lemma \ref{FormuleInverse}.
Hence (\ref{SpectralEstimat2}) can be re-written
\begin{equation}\label{SpectralEstimat3}
  c_{-}(\bu,\bu)_{\mS}\leq (\mR\mB^{-1}\mR^{*}\mA\bu,\bu)_{\mS}
  \leq c_{+}(\bu,\bu)_{\mS}\quad \forall\bu\in\mH.
\end{equation}
To conclude the proof there only remains to observe that,
since $\mR\mB^{-1}\mR^{*} = \mS^{-1}$, then $\mR\mB^{-1}\mR^{*}\mA$ is
self-adjoint with respect to both $(\;,\;)_{\mS}$ and  $(\;,\;)_{\mA}$. 
As a consequence, Lemma \ref{Minimax} combined with \eqref{SpectralEstimat3}
implies \eqref{SpectralEstimate}. \hfill $\Box$

%
%


\end{document}